\documentclass[11pt,a4paper]{article}
\usepackage{t1enc}
\usepackage[latin1]{inputenc}
\usepackage[english]{babel}
\usepackage{harvard}
\newcommand{\beq}{\begin{equation}}
\newcommand{\eeq}{\end{equation}}
\newcommand{\beqarr}{\begin{eqnarray}}
\newcommand{\eeqarr}{\end{eqnarray}}
\newcommand{\beqa}{\begin{eqnarray*}}
\newcommand{\eeqa}{\end{eqnarray*}}
\begin{document}
\thispagestyle{empty}
\begin{center}
\subsection*{The changing concept of matter in H. Weyl's thought,\\ 1918 -- 1930 }
{\em Erhard Scholz, Wuppertal}  \\[5pt]
\end{center}
\vspace{5mm}
\begin{abstract}

During the ``long decade'' of transformation of mathematical physics between 1915 and 1930, H. Weyl interacted  with  physics in two highly productive phases and contributed to it, among others, by   his  widely read book on {\em Space - Time - Matter (Raum Zeit - Materie)} (1918 -- 1923) and on {\em  Group Theory and  Quantum Mechanics (Gruppentheorie und Quantenmechanik)} (1928 - 1931). In this time Weyl's understanding of the constitution of matter and its mathematical description changed considerably. At the beginning of the period he started from a "dynamistic" and geometrical conception of matter, following and extending the Mie-Hilbert approch, which he gave up during the year 1920. After transitional experiments with a singularity  (and in this sense topological) approach  in 1921/22, he developed an open perspective of what he called an ``agency theory'' of  matter. The idea for it was formulated already before the advent of the ``new'' quantum mechanics in 1925/26. It turned out to be well suited to be taken over to the quantum view as a kind of ``heritage'' from the  first half of the decade.
 At the end of the period, Weyl completely renounced his earlier belief in  the possibility to ``construct matter'' from a  geometrically unified field theory. He now posed the 
possibility of a  geometrization of the mathematical  forms underlying the rising quantum physical description of matter as a completely open problem for future research.

\end{abstract}
\subsubsection*{Introduction}
It may appear a  strange question to ask for the changing views   of a mathematician  on the concept of matter. Why not pose it for a  natural scientist or a philosopher?  But Hermann Weyl was, as we  know,  a bit of  all of them. His views on mathematics and their foundations made it impossible to separate mathematics and its ``meaning'' from broader contexts of its use as a conceptual form and as a symbolic tool for the understanding of nature (or at least some aspects of it).  

During the ``long decade'' of transformation of mathematical physics, as we may call the time between 1915 and 1930, with the rise of the general theory of relativity (GRT) and the origin of the ``new'' quantum mechanics (QM), H. Weyl interacted  with  physics in two highly productive phases and contributed to the development of both, theoretical physics and the mathematical concepts and methods in it.  The first phase lasted from  1916  to  1923 and had as main outspring  his  widely read book on {\em Space - Time - Matter (Raum Zeit - Materie)} \cite{Weyl:RZM}, which we will also refer to by RZM.  In this period the book had five successive editions with considerable  extensions and/or alterations well documenting the shifts in the understanding of the subject by its author. Some of these changes were of a more technical nature for general relativity or the mathematics involved, others were of more basic nature,  including in particular the changing characterization and mathematical description of matter.  In the middle of the 1920s Weyl worked on the representation theory of Lie groups (1924/1925) and wrote a book on the philosophy of mathematics and the natural sciences (1925/1926), before he started to contribute  actively to the  rising quantum mechanics, culminating in his second book about mathematical physics, {\em Group Theory and Quantum Mechanics} \cite{Weyl:GQM}. 

The growing awareness of the unreducible and far-reaching role of quantum properties had already contributed to considerable shifts of Weyl's concept of matter during the first phase of involvement  in physics.  In the second half of the ``long decade'' his views were deeply transformed by the rise of quantum physics. 
 This transformation was, of course, much more than  a personal experience. It reflected the experience of the whole community of  researchers in basic physics of the time, although seen from a specific Weylian perspective. As such, it may be  illuminating for a historical and philosophical understanding of  the transformation of  the concept of matter, brought about by the   tension resulting from the unfinished  ``double revolution'' of GRT and QM during the 1920s.\footnote{A. Pais'   description of the change of matter concepts by the rising quantum theory as   ``the end of the game of pebbles''  \cite[324]{Pais:Inward} fits already well to  this shift, although Pais used  it as a header for the rise of second quantized fields starting in the late 1920s.} 
 In spite of the drastic difference between Weyl's concept of matter at the beginning of the period and at the end of it,  we easily perceive  a  common thread linking both ends. This common underlying feature is  a dynamistic  view of matter. This characterization has to be understood in a  general, philosophico-conceptual sense which {\em may} be  related, but {\em need not} be, to the electrodynamical picture of matter  which gave a new thrive to dynamism among physicists  and mathematicians of the early 20th century. 

In the history and philosophy of physics, the dynamistic view of matter  in the early 20th century is  often restricted to the exclusively electromagnetic approach. Such a restriction  shadows off the  intricate link to the quantum theoretical phase, which played a role for some of the protagonists of the period. Of course, also  Weyl started  from  Mie's electrodynamical theory of matter when he first looked for an adequate ``modern''  mathematical expression of such a dynamistic view.  From this basis  he developed his  program of a geometrically unified  field theory in the first phase of his involvement in mathematical physics.\footnote{\label{UFT}For Weyl's first phase of involvement in mathematical physics compare \cite{Skuli:Diss,Scholz:DMV}, for  broader views on unified field theories see \cite{Vizgin:UFT,Goldstein/Ritter:UFT,Cao:FTs,Goenner:UFT}. }
 The impact of quantum physics replaced classical field pictures by  quantum stochastical descriptions of  the ``agency nature'' of matter, as Weyl liked to call it.   
At the end of the period discussed here, he completely renounced his earlier belief in  the possibility to ``construct matter'' from a  geometrically unified field theory. He now posed the 
possibility of a  geometrization of the mathematical  forms underlying the rising quantum physical description of matter as a completely open problem for future research.

In this article I present a kind of ``logitudinal section'' through the long decade, oberved along the trajectory of a  single person, who was partially a contributor and partially a well informed observer of the  development.\footnote{A  complementary  view at several ``transversal'' sections (in time) with a broad evaluation of authors and approaches is presented in \cite{Goldstein/Ritter:UFT}.}
 We start with Weyl's turn towards Mie's theory of matter, his own contribution to it, and his rather early distachment, which was related to the influences of early quantum mechnics without being a  necessary conclusion from it. After a short phase of relaxation of classical explanations of matter, by a combination of metrical and topological aspects (matter characterized by singularities in  space-time), Weyl developed an open perspective of what he called an ``agency theory'' of matter. The idea for it was formulated already before the advent of the ``new'' quantum mechanics in 1925/26. It turned out to be well suited to be taken over to the quantum view as a kind of ``heritage'' from the  first half of the decade.

\subsubsection*{Adherence to Mie's dynamistic approach to matter}
After Weyl came back to neutral Switzerland from his war duty in the German army in May 1916, he started a completely new phase of his research, which was imbued by a  longing  for a sounder basis of knowledge.\footnote{See \cite[64ff.]{Skuli:Diss}, \cite{Schappacher:Politisches}.}
 For him,  this meant to   work  in a broad and  interconnected  set of   fields comprising the foundations of analysis, differential geometry,  general relativity, unified field theory   and the basic structures of matter.  Only if we take this broad range of intellectual activities into account, we can get an adequate sense of Weyl's conceptual and theoretical moves inside the single fields.  Let us have a look at some  points of such interconnections:
\begin{itemize}
\item[-] In the foundations of mathematics
our author shifted from his own construc\-tive-arithmetical approach for a characterization of the concept of continuum \cite{Weyl:Kontinuum} to a kind of Brouwerian intutionism \cite{Weyl:Krise}. For a while, he believed Brouwer's approach to possess an intimate connection to his ideas in purely infinitesimal geometry.   Weyl  could  well characterize ``purely infinitesimal'' structures  on the level of differential geometry by  his  generalization of a Riemannian metric by  combining  a conformal structure  with a {\em length connection} $\varphi = \sum_i \varphi_i dx^i$.  It was, however, much more difficult to give them a  mathematical meaning on 
  the foundational (and topological) level. Here a  precise  conceptual characterization was lacking. Weyl was well aware of this deficiency which  contributed to  tensions and shifts inside his  foundational  contributions. For a while,  Brouwer's ``revolutionary'' approach to the continuum  (as Weyl called it in 1920) appeared  him to offer a promising road.\footnote{Compare \cite[121ff.]{Hesseling:Gnomes}, \cite{Scholz:Continuum}.} 
\item[-]  For some years,  Weyl considered his gauge geometrical generalization of the Riemannian metric as the proper approach for a unified field theory of gravitation and electromagnetism and, moreover, a field theory of matter based upon it.\footnote{See footnote \ref{UFT}.}
\item[-] Rising doubts with respect to  the physical feasibility of this immediate physical interpretation of  gauge geometry contributed to a turn towards   a more basic 
philosophico-conceptual analysis of the principles of congruence geometry in Weyl's {\em mathematical analysis of the problem of space}.\footnote{\cite{Scholz:SpaceProblem}}
\item[-] The necessity, or at least usefulness, to accept classical logical prinicples (excluded middle) in the proof of the main theorem of the analysis of the space problem contributed to rethink his radical position in the foundations of mathematics.
\end{itemize}
  
In the second point indicated above, the field theoretic approach to matter constitution, Weyl was deeply influenced by  Mie's electromagnetic theory of  matter which he got to know through Hilbert's modification during the autumn  1915.\footnote{Compare \cite{Corry:Hilbert_Red,Corry:Hilbert/Mie,Kohl:Mie,Sauer:Hilbert,Vizgin:UFT}. } 
Hilbert   attempted to arrive at a  kind of mathematical synthesis  of Mie's and Einstein's ideas  on electromagnetism (Mie) and gravitation (Einstein). He indicated how to find a  united Hamiltonian for gravitation  and electromagnetism in a generally covariant setting.\footnote{For a  discussion of Hilbert's research program building upon and extending Mie's field theoretic matter theory see  \cite{Sauer:Hilbert}, for a  critical evaluation of Hilbert's relation to Einstein's theory of general relativity  \cite{Corry/Stachel/Renn,Renn/Stachel:Hilbert}.}
He was convinced, that in such a classical united field theory the riddles of the grainy structure of matter should be solvable. H. Weyl and F. Klein were not convinced that Hilbert's  attempted ``synthesis'' of Mie and Einstein was acceptable as a physical theory. They argued for  a broader understanding of Hilbert's approach.  E. Noether's mathematical analysis of Hilbert's invariance conjectures (later ``Noether theorems'') contributed an essential mathematical stepping stone for it.\footnote{See \cite{Rowe:Noether,Brading:Noether_Weyl}.}

Although a central point of the scepticism resulted from  the unclear role of energy conservation in Hilbert's approach, Weyl  was, moreover,  not  convinced that Hilbert's approach was able to lead to a unification of gravitation and electromagnetism, in which  matter structures were better derivable  than in Mie's original version. In the first edition of his book he therefore discussed  a field theoretic matter concept essentially as  in   Mie's original  purely electromagnetic approach 
 \cite[\S 25]{Weyl:RZM}. Hilbert's generalization was  only  mentioned in passing, in the section   which  treated the   modification of the   Hamiltonian principle for electromagnetism by gravitation \cite[\S 32]{Weyl:RZM}. On the other hand, Weyl  presented Mie's approach in such a convinced rhetoric form that the reader might easily get the impression that  Mie's research goal had  already nearly been achieved. The desired  result (derivation of a ``grainy'' structure from field laws) seemed close to  sure. After a comparison of  Mie's theory  with Maxwell-Lorentz's, Weyl stated:
\begin{quote}
The theory of Maxwell and Lorentz cannot hold for the interior of the electron; therefore, from the point of view of the ordinary theory of the electrons we must treat the electrons as something given {\em a priori}, as a foreing body to the field. A  more general theory of  electrodynamics has been proposed by {\em Mie}, by which it seems possible to derive the matter from the field
 \ldots .\cite[165]{Weyl:RZM}
\end{quote}
This formulation was kept unchanged by Weyl in the next three editions.\footnote{It thus appears verbally unchanged  in the third edition on which H.L. Brose's English translation is based \cite[206]{Weyl:RZMEnglish}. Here, as in other cases, our English quotes from RZM  are following  Brose's  translation, where available.  } He only  changed it during the last revision for the fifth edition (1923). Then  he  clearly expressed the open status of Mie's attempt and presented it, in an essentially didactical approach, as nothing but    an {\em example} of  a physical theory ``which agrees completely with the recent ideas about matter''   \cite[$^5$1923, 210]{Weyl:RZM}.

Mie's proposal fitted beautifully to Einstein's detection of the energy-mass equivalence of special relativity and seemed to extend it. In a passage commenting the  equivalence  $ E = mc ^2$, Weyl argued:
\begin{quote}
We have thus attained a new, purely dynamical view of matter (footnote: Even Kant in his ``Metaphysische Anfangsgr\"unde der Naturwissenschaft'' teaches the doctrine that  matter fills space not by its mere existence but in virtue of the repulsive forces of all its parts.)
Just as the theory of  relativity theory has taught us to reject the belief that we can recognize one and the same point in space  at different times, {\em so now we see that there is no longer a meaning in speaking of the  same position of matter  at different times}. \cite[162]{Weyl:RZM}, \cite[202]{Weyl:RZMEnglish}
\end{quote}
Already here, in the context of special relativity, he  described an electron as a kind of "energy knot'' which ``propagates through empty space  like a water wave across the sea'', and which could  no longer be considered as element of some  self-identical substance. Then, of course, there arose the  problem to  understand both,  this kind of propagation of energy, and the stability of the ``energy knot''. Weyl stated   the new challenge of (special) relativity to field theory, which arose  from a dynamical understanding of matter/energy:
\begin{quote}
The theory of fields has to explain why the field is granular in structure and why these energy-knots preserve themselves permanently from energy  and momentum in their passage to and fro (\ldots); therein lies the {\em problem of matter}.
\cite[162, emphasis in original]{Weyl:RZM} \cite[203]{Weyl:RZMEnglish}
\end{quote}
Like the dynamists of the early 19th century, Weyl now insisted that atoms could not be considered as  invariant  fundamental  constituents of matter:
\begin{quote}
{\em Atoms and electrons are not}, of course, {\em ultimate invariable elements}, which  natural forces seize from without, pushing them hither and thither, but they are themselves distributed continuously und subject to minute changes of a fluid character in their smallest pieces. It is not the field that requires matter as its carrier in order to be able to exist itself, but {\em matter} is, on the contrary, {\em an offspring of the field}.
 (ibid.).\footnote{Translation slightly adapted, E.S.}
\end{quote}
It seems worthwhile to remark that these general passages on the dynamistic outlook on the problem of matter  were {\em not changed} by Weyl until (and including) the fifth edition of his book in 1923. On the other hand,  the special role attributed to Mie's theory, or to his own unified field theoretic approach, underwent considerable changes during the following years. But in spite of all his enthusiasm for the new role of field theory in  the understanding of matter, Weyl indicated already in 1918 after his presentation of Mie's theory, that something new was rising at the (epistemic) horizon, which might have  unforeseen consequences in the future. He compared the actual status of field physics with the seemingly all-embracing character of Newtonian mass-point dynamics in the Laplace program at the turn to the 19th century and warned:
\begin{quote}
Physics, this time as a physics of fields, is again pursuing the object of reducing the totality  of natural phenomena to {\em a single physical law:}  it was believed that  this goal was almost within reach when the mechanical physics of mass points, founded upon Newton's  Prinicipia, was celebrating its  triumphs.  But also today, provision is taken that our trees do not grow up to  the sky. 
 We do not yet know  whether the state quantities underlying Mie's theory suffice for a characterization of matter,  whether it is in fact   purely ``electrical'' in nature. Above all, the dark cloud of all those appearances that we are  provisionally seeking to  deal with by the quantum of action throws its shadow upon the land of of physical knowledge,  threatening no one  knows what  new revolution.
\cite[170]{Weyl:RZM}\footnote{Unchanged in all editions, last one in \cite[$^5$1923, 216]{Weyl:RZM}. Translation from \cite[212]{Weyl:RZMEnglish} slightly adapted by E.S.}
\end{quote}

\subsubsection*{A geometrical extension of  Mie's theory of matter }
A few months after his book manuscript   was finished, Weyl developed his concept of a generalizated {\em Weylian metric} on a differentiable manifold.   In technical terms his metric was  given,  and still can be characterized,  by an equivalence class of pairs, $[(g, \varphi)]$, consisting of a (semi)Riemannian metrics $g = \sum_{i,j} g_{i j} dx_i dx_j$ and  a differential form (``length connection'') $\varphi = \sum_i \varphi_i dx_i$, up to equivalence by conformal factors in the Riemannian component of the  metric and ``gauge transformation'' of the length connection form.\footnote{\cite{Varadarajan:connections,Vizgin:UFT}}  This generalization allowed a seemingly natural interpretation of the potential of the electromagnetic field by the length connection and thus a metrical unification of the main physical fields known at  the time, gravitation ($g$) and electromagnetism ($\varphi$). Weyl  considered this structure as an important step forward for the Mie program of a dynamical characterization of  matter. He  published about it in several articles, starting in 1918.   In the following  year he   included the approach into the third edition of his book \cite[$^3$1919]{Weyl:RZM}.

The first edition had ended with a section on cosmology,  ``Considerations of the world as a whole'' \cite[\S 33]{Weyl:RZM}. In the third edition two new sections were added, one on ``the world metrics as the origin of the electromagnetic phenomena'',  containing an introduction to Weyl's unified field theory, and one on ``matter, mechanics and the presumable (mutma\ss{}liches)  law of the world'', in which Weyl's extension of the Mie program was sketched.  Like in the  first edition, Hilbert's extension of Mie's program was   only   indirectly mentioned in the section on the combined Hamiltonian principle of electromagnetism and gravitation. On the other hand,  the last section culminated in Weyl's own attempt to overbid both Mie and Hilbert by a derivation of the discrete ``granular'' matter structures from  his gauge invariant action principle.  In his lecture course on mathematics and the knowledge of nature  of the winter semester 1919/20 , Hilbert countered by an  acid remark  that such a perspective would lead to a kind of  ``Hegelian physics'', in which the ``whole world process would not go becond the limited content of a finite thought'' \cite[100]{Hilbert:Rowe}. He did not explain, though, why this kind of analysis should not  apply  to his own program just as well .

In 1919 Weyl was at the high point of enthusiasm for his new theory. The new section in the third edition of his book started with a rhetoric trumpet-blast:
\begin{quote}We  rise to a final synthesis \ldots \cite[$^3$1919, 242]{Weyl:RZM}\footnote{``Wir erheben uns zu einer letzten Synthese.'' Brose's translation reduced the kick of enthusiasm considerably: ``We now aim at a final synthesis'' \cite[282]{Weyl:RZMEnglish}. Weyl did not weaken the rhetoric until and including the fifth edition, although he slightly revised its wording by adding a ``nun  (now)'' \cite[$^5$1923]{Weyl:RZM}.}
\end{quote}
Part of his enthusiasm resulted apparently from the realization that gauge invariance with respect to the change of the length gauge led to a new invariance principle, which in Weyl's semantics of the approach could only be the invariance of electrical charge \cite[$^3$1919]{Weyl:RZM} \cite[293]{Weyl:RZMEnglish}. That was, of course, a great achievement of lasting importance, even if the specific version of gauge invariance had later to be given up.\footnote{See \cite{Vizgin:UFT,Brading:Noether_Weyl} and for a detailed, historically oriented discussion of the underlying  mathematics \cite{Varadarajan:connections}.} But Weyl hoped for more. He expected that on the one hand the cosmological modification of the Einstein equation should be a natural result from his gauge geometry. On the other hand the stable solutions of the equations for the ``problem of matter'', satisfying adequate regularity conditions should lead to a  discrete set of solutions depending on some parameter $\beta $. This expectation  had a (formal)  similarity to a set of ``discrete eigenvalues'' of an operator, although here the operator was not linear.  

The problem was, in fact,  characterized by a non-linear differential equation of great   complexity. Even Weyl guessed that the available tools of analysis would probably  neither suffice for a proof of their existence, nor for an approximative calculation \cite[$^3$1919, 260]{Weyl:RZM}. This remark made the epistemic status of Weyl's ``discrete solutions''  highly problematic. It turned them rather into a  symbol for a  natural philosophical speculation than into an object for research in  mathematical physics. Weyl continued the discussion by a beautiful  remark.
\begin{quote}
The corspusculae which correspond to the possible eigenvalues had to coexist in the same world besides  each other or in another,   mutually enforcing on another  subtle modifications of their intrinsic structure; strange  consequences seem here to arise  for the organization of the universe;    perhaps they  may make comprehensible its stillness in the large and  unrest in the small. \cite[$^3$1919, 261]{Weyl:RZM}
\end{quote}

When Weyl wrote these lines, he was at the peak of his belief in a strong unification program of forces and matter, which could be constructed on purely geometrical grounds.  In the last long passage of the new added sections we find a discussion of how he now saw  the relationship between geometry and physics in the light of his recent findings.
\begin{quote}
We have realized that physics and geometry coincide with each other and that the world metrics is one, and even the only one, physical reality. Thus, in the final consequence, this physical reality appears   as nothing but a pure form; geometry has not been physicalized but physics has been geometrized (nicht die Geometrie ist zur Physik, sondern die Physik ist zur Geometrie geworden). \cite[$^3$1919]{Weyl:RZM}
\end{quote}
H. Weyl was now at the apogee of the belief in a  strong unification which was both,  deeply reductionist and highly idealistic. In his eyes, physics seemed  to be transformed to a purely formal status and was absorbed by geometry. Matter had  seemingly become an epiphenomenon of the ``world metrics'' which started to acquire  a slightly mystical flavour.  In the mind of our protagonist, the physicalizing tendency of geometry among leading protagonists of the 19th century, including researchers like C.F. Gauss, N.I. Lobachevsky, and B. Riemann, appeared turned upside down  ---  even though  in this extreme form for only  a year or two. 

As we will see in  a moment, this conviction did not hold for long. Already in the fourth edition, this extremely reductionist  passage was cancelled by its author. Now the book ended with another, less reductive passage on the unifying power of the mind  and  an \'eloge of the  ``chords from that harmony of the spheres of which Pythagoras and Kepler once dreamed'' \cite[312]{Weyl:RZMEnglish}. Weyl did not hide that he had  changed his mind; in a separate article written for the physical community shortly after the revisions for the fourth edition, he explained frankly:
\begin{quote}
From the first edition of RZM to the third one I took the position of ( \ldots ) [ a purely field theoretic characterization of matter, E.S.], as I was charmed by the beauty and unity of pure field theory; in the fourth edition, however, I lost confidence in the field theory of matter by striking reasons and changed to the second point of view [of a primacy of  matter,  irreducible  to interaction fields, E.S.] . \cite[242]{Weyl:47}
\end{quote}
Let us therefore have a look at the 
  ``striking resons'' for this ontological shift at the beginning of the 1920s.

\subsubsection*{From speculations on the ``the causal and the statistical view of physics'' to a break with Mie's theory of matter }
In the year 1919 Weyl gave a talk to the Swiss {\em Naturforschende Gesellschaft} on the relationship between the causal to the statistical view of physics which was published a year later \cite{Weyl:Kaus}. This  paper has been strongly critcized by  P. Forman in his otherwise very stimulating article on Weimar culture and its influence on the discourse among physicists \cite{Forman:I} as a document of an ``antirational'' kind of ``conversion to acausality''. We need not take up here again the broader debate on the question how ``antirational'' the move was and what kind of ``acausality'' was at stake here.\footnote{See the illuminating comments and critique of Forman's original presentation in \cite{Hendry,Skuli:Diss,Stoeltzner:Indeterminism} and also the modifications in \cite{Forman:II}.}
It  may suffice to add that for Weyl the topic of his talk contained a challenging combination of questions in the conceptual foundations of contemporary physics, including the rising ``clouds'' of quantum phenomena, with the question of  how modern natural science can be made compatible with metaphysical considerations of the existential experience of the openness of evolving life processes  and of the freedom of personal actions. A central topic of this talk was the directedness and irreversibility of time, which  appeared Weyl to be linked to some process  level of irreducibly statistical  nature.\footnote{Weyl hinted at the possibility that the classical mechanical discussion on ergodicity had to be revised the light of ``some mysterious discontinuity'' introduced recently by quantum theory \cite[118]{Weyl:Kaus}.} 

The talk took place about the time of his turn  to Brouwer's intuitionism.  In {\em this respect} we  could even speak of a  kind of  ``conversion''.\footnote{Compare, e.g., \cite[127]{Hesseling:Gnomes}.} Weyl speculated that perhaps  Brouwer's approach could lead to a solution of several fundamental  problems at one strike. In mathematics he hoped  for an answer to the foundational question of the concept of the mathematical continuum and  for a philosophically and mathematically sound characterization of the topological ``substrate'' of purely infinitesimal geometry; in physics he expected  a break with the rigidity of the  causality structure in classical mechanics (``Gesetzesphysik'')  and an access to understand the irreducible directness of time. He expected that satisfying answers to all these questions might have some intimate link to  the open process of ``becoming'' inherent in Brouwer's choice sequence for the characterization of the intuitionistic continuum:
\begin{quote}
Finally and foremost, it is inbuilt into the essence of the continuum that it cannot be treated as a rigid being, but rather only as something what is continuously evolving in an infinite, inward bound  process of becoming. \cite[121]{Weyl:Kaus}
\end{quote} 
In this speculative thought, Weyl hoped to find a common thread binding together  the foundations of  analysis, ``purely infinitesimal''  geometry,  the  directedness of time flow in the physical world,   its determinative openness, and  a conjectured irreducibly  stochastical nature of physical laws, which would break with the classical kind of lawfulness (``causality'' in the language of the time).

In September 1920, during the discussions of  the Bad Nauheim meeting of the German {\em Naturforscher Versammlung} and from a draft manuscript of Pauli's contribution on relativity to the {\em Enzyklop\"adie der Wissenschaften} Weyl got to know content and reason of  Pauli's critical evaluation of his  modified version of  the Mie theory of matter. This conjunction of detailed scientific criticism, coming from a personally close, young expert in the field, with his own most recent  conceptual and metaphysical speculations, undermining the classically deterministic  field structures anyhow, led Weyl to give up the belief in his program of  a geometrically unified field theoretical derivation of matter structures. At the end of the year, in a letter to Felix Klein, in which he reported on his recent advances on mathematical and physical questions (included or not into the just finished fourth edition of RZM), he reported among others:
\begin{quote}
Finally I thoroughly distached myself from Mie's theory and came to a different position with respect to the problem of matter.  I no longer accept field physics  as the key to reality. The field, the ether, appears to me only as a {\em transmitter} of effects, which is completely feeble by itself; while matter is a reality lying beyond the field and causing its states \ldots . \cite{WeylanKlein:12/1920}\footnote{``Endlich habe ich mich gr\"undlich von der Mie'schen Theorie losgemacht und bin zu einer anderen Stellung zum Problem der Materie gelangt. Die Feldphysik erscheint  mir keineswegs mehr als der Schl\"ussel zu der Wirklichkeit; sondern das Feld, der \"Ather ist mir nur noch der in sich selbst v\"ollig kraftlose {\em \"Ubermittler} der Wirkungen, die Materie aber eine jenseits des Feldes liegende und dessen Zus\"ande verursachende Realit\"at. Mit dem ``Weltgesetz'' (Hamiltonsches Prinzip), das die Wirkungs\"ubertragung im \"Ather regelt, w\"are noch gar wenig f\"ur das Verst\"andnis aller Naturerscheinungen gewonnen.'' (ibid. emphasis in original); compare \cite{Skuli:Diss}.}
\end{quote}
Similar phrases are to be found close to the end of  the fourth edition of RZM. Here the last section, containg Weyl's version of Mie's theory was no longer announced under the emphatic title ``Matter, mechanics, and presumable  world law'' as in the third edition. The discussion  was now only presented  as a ``development of the simplest principle of action \ldots''. 

It contained a short discussion of some consequences of Weyl's  gauge invariant quadratic  action $S^2 \sqrt{|det g|}$ for the Hamiltonian of a combined theory of gravitation and electromagnetism, with $S$ the scalar curvature of Weyl geometry. Now he commented that this action is only the ``simplest assumption for calculation'', for which the author no longer wanted to  ``insist that it is realised in nature'' \cite[295]{Weyl:RZMEnglish}. For anybody who continued to read the book until the end, Weyl  made  clear that he now conjectured  a close interrelation between the directedness of time flow with quantum jumps as seen  in the Bohr model of the atom. That was no longer compatible with the classical structures of time invertible determinism:
\begin{quote}
We must here state in unmistakable language that physics at its present stage can in no wise be regarded as lending support to the belief that there is a causality of physical nature  which is founded on rigorously exact laws. The extended field, ``ether'' is merely the {\em transmitter} of effects and is, of itself, powerless; it plays a part that is in no wise different from that which space with its  rigid Euclidean metrical structure plays according to the old view; but now the rigid motionless character has become transformed into one which gently yields and adapts itself. \ldots  \cite[311, emphasis in original]{Weyl:RZMEnglish}
\end{quote}
Now the old duality of field (``ether'') and matter was back again  for our protagonist. This brought him closer to the perception of the problem by the majority of physicists working on the structure of matter and indicated a growing distance to the views held by A. Einstein. 

\subsubsection*{A short-lived singularity theory making place to  ``agency structures''   of matter}
As Weyl came from  as strong field theoretic paradigm, it was natural for him in the years 1920/21 to characterize  matter   by its formal relationship to the interaction field(s).\footnote{For Weyl, his unification of the interaction  still seemed valid, after his distachement from Mie's matter theory, and thus the singular form ``field'' would apply;  scepticists  with respect to his unification  could easily reread his remarks by turning to the plural ``fields''. }  
Thus in the fourth edition of RZM Weyl  stated  his new viewpoint clearly:
\begin{quote}
Contrary to Mie's view, {\em matter} now appears {\em as a real singularity of the field}. \cite[262, emphasis in original]{Weyl:RZMEnglish}\footnote{In the German original \cite[$^4$1921, 238]{Weyl:RZM}, no longer in the fifth edition.}
\end{quote}
But then, matter had somehow to be located in a determinative boundary structure of the field and the old question of the  structures of matter, which in the classical mechanistic view had been  given by the assumptions of its atomic constitution and the hypothetical extension of mechanical laws to the  atomic level,  was again open. After the experience of  dynamistic hopes during his period of adherence to the Mie theory, and in the light of recent  modifications coming from experimental knowledge in microphysics, Weyl  came to the conclusion:
\begin{quote}
If matter is to be regarded as a boundary singularity of the field, our field-equations make assertions only about the {\em possible states of the field} and {\em not about the conditioning of the states of the field by the matter}. This gap is provisionally filled by the {\em quantum theory} in a manner of which the underlying principles are still completely ununderstood. \cite[303, emphasis in original]{Weyl:RZMEnglish}\footnote{The translation of the last phrase has been slightly changed to adapt it closer  to the  German original \cite[$^4$1921, 276]{Weyl:RZM} than in Brose's translation. This passage is  no longer contained in the fifth edition.}
\end{quote}

Now the task to understand matter mathematically could be approached from different viewpoints. One was topological in nature. General relativity offered the opportunity to consider a differential topological manifold with boundaries, in the interior of which the fields are regular, while they are singular on the boundaries and  diverge in respective limiting processes. In an article written for {\em Annalen der Physik} shortly after the publication of the fourth edition of RZM, Weyl explained his new viewpoint more in detail \cite{Weyl:47}. He argued in two directions. Coming from the point of view of special relativity and Minkowski space, the generalization for GRT consisted not only in a deformation of the metric but could also comprise a topological modification of the underlying manifold. Weyl argued that from the space-time manifold with a combined electromagnetic and gravitational field, the subsets on which the fields obtain singular values should be omitted.
\begin{quote}
In the general theory of relativity the world can possess arbitrary (...)  connectedness: nothing excludes the assumption that in its Analysis-Situs properties it behaves like a four-dimensional Euclidean continuum, from which different tubes of infinite length in one dimension are  off. \cite[252f.]{Weyl:47}
\end{quote}
  If the general relativistic point of view was considered as the more realistic one, it even appeared as more natural to turn the view round. One  would then   have to argue in terms of pasting rather than of cutting:
 \begin{quote}
The simply connected continuum from which we construct the domain of the field by cutting off the  tubes is nothing but a mathematical fiction, although the metrical relations persisting in the field strongly propose the extension of the real space by  addition of such fictitious improper  (erdichteter uneigentlicher) regions  corresponding to the single matter particles. (ibid.)
\end{quote}
For Weyl, this change of the mathematical construction of space went in hand with a change of the understanding of the relationship between space-time and matter:
\begin{quote}
According to [this] perception, {\em matter itself is nothing spatial (extensive) at all, although it is inserted in a certain spatial neighbourhood}. \cite[254, emphasis in original]{Weyl:47}
\end{quote}
He must have liked this idea. One of the Fichtean motifs on the ``construction'' of matter and space from ``forces'', which had impressed Weyl already at the time of his turn towards ``purely infinitesimal'' geometry,  ancquired here a new face and persisted  in a modified form.\footnote{Compare \cite{Scholz:ICM,Scholz:Cret}.}

On the other hand, there was a physical approach to the problem of matter. In addition to the proper laws of the field(s)  one  had to ``study the laws according to which matter excites the field actions''.  For Weyl, matter was now turning into an irreducible originator of dynamical excitation of the interaction  field(s) and was itself guided by the latter in its own  spatio-temporal dynamics. He insisted that it could neither  be understood as a  ``substance'' in the  sense of   traditional natural philosophy, nor could it be derived from the ``field'' as in the Mie version of  dynamicist matter explanation. Weyl preferred to give a description and a term of his own to characterize matter as an {\em agency (Agens)}. This word was  uncommon in the German  language, and,  to my knowledge,  even not used in the earlier discourses on natural philosophy. Weyl   apparently shaped it on his own from the Latin gerundial form  {\em agens} for something that is acting.  In his usage of the word, an agency perception of matter was not far from the older dynamicist one of the philosophical debate of the early 19th century. As, however, the dynamistic view of matter had been linked to the electromagnetic world view and its generalizations in the Mie - Hilbert - Weyl approach by classical field theories, Weyl had good reasons to demarcate the break with this semantical field by the choice of a new word.

Different to the older field theoretical  theories, Weyl considered it as an important feature of  the agency view that it considered  matter as something which acts on spatial structures like fields, although it is not itself located inside space. Already in 1921, several years before the advent of the refined form of the quantum mechanics, Weyl stated optimistically:
\begin{quote}
In addition to  the {\em substance} and the {\em  field} perceptions  we have to add a   third 
view of matter as an {\em agency (Agens)} effecting the field states. \ldots
It makes place for the modern physics of matter, working with statistical concepts, besides the strictly functional physics of a classical field. \cite[255]{Weyl:47}
\end{quote}
For Weyl, such a shift had nothing to do with a  longing for ``acausality'', or even the adoration of it. He rather insisted that the view of  mechanical and classically field theoretic physics had reduced causality  to a purely functional mathematical relationship, while the agency perception opened a possibility to understand the causation of field states  by matter in a new and deeper way.
\begin{quote}
Here the specified direction of the passage of time: past $\rightarrow$ future, which cannot find its place in field physics, can be taken up again; in fact it is most closely related to the idea of causation. \cite[256]{Weyl:47}
\end{quote}
The characterization of causality by a deterministic and time-invertible lawlike structure as in classical mechanics appeared him as an inappropriate concept. The change from classical determination to a probabilistic one would therefore not at all contradict the concept of ``causation''. Just to the contrary, Weyl expected that it might open the path towards a more appropriate understanding of the latter. Although his most recent turn  had its origin in the short-lived singularity theory of matter, the ``agency paradigm'' of matter was kept open for  a  modification in its mathematical characterization and for a  future enrichment by an improved  understanding of  its physical properties. 

\subsubsection*{The agency concept of matter as an open research field }
The role Weyl assigned to singularities of classical fields in the fourth editionoif RZM remained itself ``singular'' in his work. It did not appear earlier and vanished, or was at least drastically reduced in importance, nearly as fast as it appeared. 
In the fifth edition of his book the section on ``further rigorous solutions of the statical problem of gravitation'', which contained the central passages on the singularity theory of the electron, from which we quoted above, was completely reorganized.
 Apparently Weyl was not satisfied with the outlook on the strong interpretation of singularities as {\em the} mathematical clue to the solution of the ``problem of matter''. In the fifth edition and in his later publications on the philosophy of nature \cite{Weyl:Materie,Weyl:PMN} we find the singularity model only in a weak sense, mentioned in passing as nothing more than an idea  illuminating the impossibility to    localize the basic agency structures of matter directly. 

In the fifth edition of RZM, Weyl no longer gave the impression that he was already in possession of a  mathematical clue to the solution of the ``problem of matter''. He now preferred  to characterize only   the terrain of investigation and discussed different approaches that had been tried, up to then.  Among these he mentioned, of course,  Mie's theory and his own generalization as important examples. But now they were   only presented as explorative theoretical models, without any claim that they might lead towards a reliable  representation of reality.

 In this discussion we find  beautiful, nearly poetic descriptions of the actual state of knowledge as an open terrain:

\begin{quote}
We only percieve the bounding embankment of the subtle, deep groove   which is dug  into the metrical face  of the world by  the trajectory of the electron; what is 
covered  by the depth,   remains hidden to us.  It may be that the whole groove is filled by a field, qualitatively equivalent to the outer one, as Mie assumed; {\em but just as well  the abyss  may be  fathomless.} Mie's perception dissolves matter into the field; the other one removes   it, so to speak, from the field. According to the latter view {\em matter is an agency determining the field, although in itself nothing spacelike, extensional, but only located in a certain spatial neighbourhood,} from which its field effects depart. \ldots \\ \cite[ $^5$1923, 286]{Weyl:RZM}
\end{quote}

Coming closer to the middle of the 1920s, Weyl left it open whether  it seemed more promising to smooth off the field   for a mathematical representation of the basic constituents of matter (like in the Mie approach),  to excise it (like in the singularity approach of 1921), or to find any other characterization which might take the statistical nature of quantum descriptions better intor account than the other ones: 
\begin{quote}
Our description of the field surrounding an electron is a first, stuttering formulation of such laws. Here lies the working field for modern physics of matter, to which  belong, above all, the facts and riddles of the quantum of action  (\ldots) As far as we can judge today, the lawfulness acording to which matter induces effects can   be described in statistical terms only, \ldots .'' (RZM  $^5$1923, 286f.)
\end{quote}
Independent of these open problems for an adequate mathematical characterization of matter, it  now appeared  clear to him  that  matter, rather than the field had to be given primacy for all  experimental purposes or any practical exchange with nature. 
\begin{quote}
Our willful actions have,  primarily, always to grapple on matter, only thus we can change the field. In fact, we then need two kinds of laws for the explanation of natural phenomena:
1. {\em field laws} (\ldots), 2. {\em laws regulating the excitation of the field by matter. } \ldots . (ibid.) \end{quote}
Finally  causality came back, for Weyl,   to be a relation which enabled human beings to influence  the course of  natural processes by a willful modification of material constellations in the world.
As  he had come to the insight that physical  knowledge of the basic matter structures was still highly restricted,  he  cancelled those  passages of the final sections of earlier editions of RZM, which appeared much too enthusiastic from  his recently acquired view. That did not exclude poetic allusions. The fifth edition, the last one revised or extending by himself, ended with a  passage which was both, sober and prophetic:
\begin{quote}
We were unable to pursue our analysis of space and time without studying matter in detail. Here, however, we are still confronting riddles the solution of which is not to be expected from field physics. In the darkness still surrounding the problem of matter, quantum theory may perhaps be the first twinkling of light. \cite[$^5$1923, 317]{Weyl:RZM}
\end{quote}
Weyl had entered the first phase of active intervention into mathematical physics (the ``RZM-phase'', as we might call it) with a strong program of reductionist unification; at the end of it, he clearly saw the necessity to distinguish ontologically and mathematically between interaction fields and matter. While for the first class the classical field theories could be considered as very successful, the problem of matter had turned back again into a riddle. 

\subsubsection*{A view back in 1930 }
Only two years after these lines were written, the  ``first twinkling of light'' was  stabilized by the establishment of quantum mechanics in the form of wave mechanics and operator theory in Hilbert spaces. The core of this development was the product of a new generation of physicists (W. Heisenberg, W. Pauli, P. Jordan, P.A.M. Dirac, E. Schr\"odinger, e.a.) who stood in close communication with outstanding figures of the earlier period (N. Bohr, M. Born, A. Sommerfeld, P. Ehrenfest, e.a.). Although Weyl was no member of this group, he  was close enough to several of the participants that he was immediately drawn into  the turn to ``the new'' quantum theory at the middle of the 1920s.  In oral and written exchange with E. Schr\"odinger, W. Pauli,  M. Born and P. Jordan he  even contributed in certain respects to  it.\footnote{For the Born and Jordan part see \cite{Scholz:IGQM}.} 
In  his lecture course in winter semester 1927/28  on {\em Group Theory and Quantum Mechanics},  he took up Schr\"odinger wave functions  and Pauli spinors (in later terminology) as new mathematical forms to represent a stochastically determining matter ``agency''. In the book arising from it \cite{Weyl:GQM} he could already integrate Dirac spinors for the characterization of a relativistic matter ``field'' of a new type. The second edition (1931) entered into the complex and irritating discussion of ``second quantization'' of these new provisional symbolic systems for the agency characterization of matter.

Knowing well about the provisional character of the quantum mechanical characterizations of matter, Weyl was deeply impressed by its successes already on the level of spectroscopy and the first steps into  the quantum chemical theory of valence bonds. An invitation to give the 1930  Rouse Ball lecture at Cambridge gave Weyl the opportunity to review the whole development of matter concepts, which had taken place during the long decade  just coming  to an end.

Even from hindsight, he still considered the attempts  of the early 1920s to geometrize ``the whole of physics'' as very comprehensible at its time, because they had tried to follow up  on Einstein's  successful geometrization of  gravitation  \cite[338]{Weyl:Rouse_Ball}. In this historizing perspective, he saw no reason to distance himself from  his own attempts of 1918. He summarized its critical reception by physicists and reviewed Eddington's approach to unification by affine connections, including  Einstein's later suppport for that program.
Comparing the latter with his own   ``metrical'' unification of 1918 he concluded that from hindsight both theory types appeared as   ``merely   geometrical dressings (geometrische Einkleidungen) rather than as proper geometrical theories of electricity''. He discussed the struggle between the metrical and affine  field theories  (i.e., Weyl 1918 versus Eddington/Einstein)  and gave the whole story a smilingly ironic turn:
\begin{quote}
\ldots there is no longer the question which of the two  theories will prevail in life, but only whether the two  have  to be buried as twin brothers in the same grave or in two different graves. \cite[343]{Weyl:Rouse_Ball}
\end{quote}

In the light of his changed view on the problem  of matter, he could find  just as little  
 arguments in favour of the more recent brands of unification attempts proposed at the end of the decade,  Einstein's distant parallelism approach  or the Kaluza-Klein approach.\footnote{For Einstein's distant parallelism see \cite{Sauer:DP}, for Kaluza and Kaluza-Klein \cite{Wuensch:2003}.} 
 Weyl completely rejected Einstein's new theory, not only by semantical reasons, but also by a mathematical one, because in his opinion it would break with the infinitesimal point of view, and warned: 
\begin{quote}
The  result [of pursuing Einstein's {\em Fernparallelismus} approach, E.S.] is to give away nearly all what has been achieved in the transition from special to general relativity. The loss is not compensated by any concrete gain.'' \cite[343]{Weyl:Rouse_Ball}
\end{quote}

Weyl percieved a  nearly complete scientific devaluation of all  unified field theories invented during  the long decade. This devaluation     resulted from  the quantum theoretical insights into matter structures, which had found first well formed mathematical representations by complex scalar or spinor fields  during the second part of the decade:
\begin{quote}
In my opinion the whole situation has changed during the last 4 or 5 years  by the detection of the matter field. All these geometrical leaps (geometrische Luftspr\"unge) have been premature, we now return to the solid ground of physical facts. \cite[343]{Weyl:Rouse_Ball}
\end{quote}
He continued to sketch the theory of spinor fields and the new understanding ot the underdetermination of phase which opened a new theoretical frame   for the  gauge prinicple. In 1929  he and V. Fock had proposed  a  revised gauge theory of electromagnetism in this context.
He insisted that the new principle of phase gauge ``has grown from experience and resumes a huge treasury of experimental facts from spectroscopy'' (ibid. 344).  That stood in marked contrast to the purely speculative principles on which all the  classical unified field theories had been built, his own one from 1918  included.
 Now he no longer expected to achieve knowledge on natural processes by geometric speculation,  but tried to anchor it in  more solid grounds, the observation of matter processes and their mathematization: 
\begin{quote}
By the new gauge invariance the {\em electromagnetic field now becomes a necessary appendix of the matter field, as it had been attached to gravitation in the old theory. } \cite[345, emphasis in original]{Weyl:Rouse_Ball}
\end{quote} 

In short, Weyl had turned from his speculative and strongly idealist approach to matter, pursued at the turn to the 1920s,   {\em to a mathematically empiristic and moderately  materialistic} one at the end of the decade. He was well aware that great   difficulties had still to be surmounted to come to grips with a quantization of the semiclassical fields (complex scalar or spinor wave functions) which had recently been invented for a provisional and partial representation of the quantum properties of matter. That gave geometry a   completely different outlook to the one in the classical field theories, although he did not want to  exclude that geometrization might  become possible some day  on a new level.   But if one wanted to continue along  this path, he was sure that  ``one had to set out in search of a geometrization of the matter field'' itself. If  one would  try to do without  an improved mathematization of the agency structures of matter themselves,  the geometrical theories would fall back to the methodological status of the unification attempts of the 1920s. He now considered  these as immature, although comprehensible first attempts, as  {\em Luftspr\"unge} (leaps into the air).  

It  may be appropriate to add that the German word ``Luftspr\"unge'' not only connotes  unrealistic first attempts, but also the joy of youthfulness. Weyl has had both, the joy of the youthful speculation to be close to a reduction of  physics to geometry,  and the   maturizing awareness that the difficult practices of experimentation and  closely related mathematical  theory production of quantum theory contained  a much  more reliable contribution towards the understanding of  the  ``agency structures'' of matter.

\subsubsection*{ }

\small
 \bibliographystyle{apsr}
  \bibliography{a_litfile}

\subsubsection*{ }
\end{document}